\documentclass[11pt, letterpaper, oneside]{myart}
\usepackage{amsmath,amssymb,amscd,amsthm,eucal, graphicx, array}
\usepackage[all]{xy} 
\usepackage[latin1]{inputenc}

\usepackage[backgroundcolor=white, bordercolor=red, linecolor=red, textsize=tiny]{todonotes} 
\usepackage{txfonts} 

\usepackage{pifont} 

\usepackage{color}
\definecolor{linkcol}{rgb}{0, 0, 0.7}
\usepackage[ linkcolor=linkcol, citecolor=linkcol, urlcolor= linkcol, colorlinks=true, hypertexnames=false]{hyperref}

\usepackage[textsize=tiny]{todonotes}

\usepackage{tikz}
\usetikzlibrary{matrix, calc,through,intersections, arrows, decorations.markings}

\renewcommand{\labelenumi}{\theenumi )} 

\numberwithin{equation}{section}
\renewcommand{\theequation}{\arabic{section}-\arabic{equation}}

\swapnumbers

\theoremstyle{plain}
\newtheorem{theorem}{Theorem}[section]
\newtheorem{lemma}[theorem]{Lemma}
\newtheorem{proposition}[theorem]{Proposition}

\theoremstyle{definition}

\newtheorem{note}[theorem]{Note}

\newcommand{\Ra}{\Rightarrow}

\newcommand{\lra}{\longrightarrow}

\newcommand{\FF}{\mathcal{F}}
\newcommand{\CC}{\mathcal{C}}

\newcommand{\Rfd}{\mathcal{R}}

\newcommand{\colim}[1][]{\ifmmode\ifinner {\operatorname{colim}_{#1}}\,\else \underset{#1}{\operatorname{colim}}\, \fi\fi}
\newcommand{\hocolim}[1][]{\ifmmode\ifinner {\operatorname{hocolim}_{#1}}\,\else \underset{#1}{\operatorname{hocolim}}\,\fi\fi}
\newcommand{\holim}[1][]{\ifmmode\ifinner {\operatorname{holim}_{#1}}\,\else \underset{#1}{\operatorname{holim}}\,\fi\fi}

\newcommand{\id}{{\ensuremath{\rm id}}}

\newcommand{\Z}{\mathbb{Z}}

\newcommand{\benu}{\begin{enumerate}}
\newcommand{\eenu}{\end{enumerate}}

\newcommand{\bit}{\begin{itemize}}
\newcommand{\eit}{\end{itemize}}

\title{Fixed points of the equivariant algebraic $K$-theory of spaces}

\author[B. Badzioch]{Bernard Badzioch}
\address[]{Department of Mathematics, University at Buffalo, SUNY, Buffalo, NY}
\email{badzioch@buffalo.edu}

\author[W. Dorabia{\l}a]{Wojciech Dorabia{\l}a}
\address[]{Department of Mathematics, Penn State Altoona, Altoona, PA }
\email{wud2@psu.edu}

\begin{document}

\date{\bf 2016.05.16}

\maketitle


\begin{abstract}
In a recent work Malkiewich and Merling proposed a definition of the equivariant $K$-theory of spaces 
for spaces equipped with an action of a finite group. We show that the fixed points of this 
spectrum admit a tom Dieck-type splitting. We also show that this splitting is compatible 
with the splitting of the equivariant suspension spectrum. The first of these results has been 
obtained independently by John Rognes.
\end{abstract}

\section{Introduction}
\label{INTRODUCTION SEC}

\footnote{Mathematics Subject Classification (2010): 19D10}
In \cite{MalkiewichMerling} Malkiewich and Merling proposed a definition of the equivariant $K$-theory of spaces 
$A_{G}(X)$ in the case when $G$ is a finite group and $X$ is a $G$-space. They also showed that 
the fixed spectrum of $A_{G}(X)$ can be described as follows. We will say that a space $Y$ is
$G$-retractive over $X$ if $Y$ is a $G$-space equipped with $G$-equivariant maps 
$r\colon Y \leftrightarrows X \colon s$ such that $s$ is a $G$-cofibration and $rs = \id_{X}$. 
Let $\Rfd^{G}(X)$  denote the category of $G$-retractive spaces  dominated  
by finite relative $G$-CW complexes
with $G$-equivariant maps over and under $X$ as morphisms. This is a Waldhausen category with weak 
equivalences and cofibrations given by $G$-homotopy equivalences and $G$-cofibrations. 
Applying  Waldhausen's $\mathcal{S}_{\bullet}$-construction to $\Rfd^{G}(X)$ we obtain a spectrum
$A^{G}(X)$.  This spectrum can be identified with the fixed 
point spectrum of $A_{G}(X)$. 

The goal of this note is to  investigate some properties of the spectrum 
$A^{G}(X)$. First, we will show that $A^{G}(X)$ admits a tom Dieck-type splitting. For a group $G$
and a subgroup $H\subseteq G$ let $NH$ denote the normalizer of $H$ in $G$ and let $WH = NH/H$.  
Denote also by  $C_{G}$  the set of all conjugacy classes $(H)$ of subgroups of $G$. 

\begin{theorem}
\label{MAIN1 THM}
Let $G$ be a finite finite group. For any $G$-space $X$ there is  a weak equivalence of spectra 
$$\tau^{A}_{X} \colon A^{G}(X) \overset{\simeq}{\lra} \prod_{(H)\in C_{G}} A(EWH\times_{WH}X^{H})$$
Moreover, this weak equivalence is natural in $X$. 
\end{theorem}

This result has been obtained independently by John Rognes.  

Next, recall that in a non-equivariant setting for a space $X$ we have the assembly map 
$a_{X}\colon Q(X_{+}) \to A(X)$ where $Q(X_{+})$ denotes the suspension spectrum of $X$.  
We will show that similarly, for a finite group $G$ and a $G$-space $X$ there is a natural map 
$a^{G}\colon Q^{G}(X_{+}) \to A^{G}(X)$ where  
 $Q^{G}(X_{+})$ is the fixed point spectrum of the equivariant suspension spectrum of $X$.
 Since $Q^{G}(X_{+})$ admits tom Dieck-type splitting a natural question is if this splitting is compatible 
 with the splitting of $A^{G}(X)$ given by Theorem \ref{MAIN1 THM}. We will show that this is the case: 
 
\begin{theorem}
\label{MAIN2 THM}
Let $G$ be a finite finite group. For any $G$-space $X$ the following diagram is commutes up to 
homotopy:
\begin{equation*}
\begin{tikzpicture}[baseline=(current bounding box.center)]
\matrix (m) 
[matrix of math nodes, row sep=3em, column sep=3em, text height=1.5ex, text depth=0.25ex]
{
Q^{G}(X_{+})&   \prod_{(H)\in C_{G}} Q(EWH\times_{WH}X^{H} _{+})  \\
A^{G}(X) &  \prod_{(H)\in C_{G}} A(EWH\times_{WH}X^{H})  \\
};
\path[->, thick, font=\scriptsize]
(m-1-1)
edge node[anchor=east] {$a^{G}$} (m-2-1)
edge node[anchor=south] {$\simeq$} node[anchor=north] {$\tau^{A}_{X}$}(m-1-2)
(m-2-1)
edge node[anchor=north] {$\simeq$} (m-2-2)
(m-1-2)
edge node[anchor=west] {$\prod a_{H}^{G}$}  (m-2-2)
; 
\end{tikzpicture}
\end{equation*}
where $a^{G}_{H}\colon Q(EWH\times_{WH}X^{H} _{+}) \to A(EWH\times_{WH}X^{H})$ is the assembly map.
\end{theorem}
 
\begin{note}
Throughout this paper we will freely use the machinery of \cite{Wal}. In particular by a Waldhausen
category we mean a category with cofibrations and weak equivalences defined in \cite[\S 1.2]{Wal}. 
Given a Waldhausen category $\CC$ by $K(\CC)$ we will denote the $K$-theory spectrum of $\CC$ obtained 
by applying Waldhausen's $\mathcal{S}_{\bullet}$-construction. 
 
\end{note}


\section{Proof of Theorem \ref{MAIN1 THM}}
\label{PROOF MAIN1 THM SEC}

Let $G$ be a finite group. Recall that $C_{G}$
denotes the set of conjugacy classes $(H)$ of subgroups  of $G$ . 
Given a subgroup $H\subseteq G$ denote by 
$\Rfd^{G}_{H}(X)$  the full subcategory of $\Rfd^{G}(X)$ whose object are all retractive $G$-spaces 
$Y$ over $X$ such that all orbits of the action of $G$ on $Y\setminus X$ are of the type $G/H$.  
This is a Waldhausen category with weak equivalences and cofibrations inherited from $\Rfd^{G}(X)$.  

\begin{proposition}
\label{AG SPLITTING PROP}
Let $G$ be a finite group, and let $X$ be a $G$-space.  There is a natural equivalence, natural in $X$:
$$\tilde{\tau}^{A}_{X} \colon A^{G}(X) \overset{\simeq}{\lra} \prod_{(H)\in C_{G}} K(\Rfd^{G}_{H}(X))$$
\end{proposition}

\begin{proof}
Consider  $C_{G}$ with an ordering:
$$C_{G} = \{(H_{1}), (H_{2}), \dots, (H_{n})\}$$
such that if $H_{j}$ has a subgroup conjugate to $H_{i}$ then $i\geq j$. 
Denote by $\Rfd^{G}_{\leq i}(X)$ the subcategory of $\Rfd^{G}(X)$ consisting of 
all retractive spaces $Y$ such that each orbit of the action of $G$ on $Y\setminus X$ is 
isomorphic to $G/H_{j}$ for some $j\leq i$. 
We claim that for any $i\leq n$ there exists a weak equivalence
$$K(\Rfd_{\leq i}^{G}(X))\overset{\simeq}{\lra} K(\Rfd_{\leq i-1}^{G}(X))\times K(\Rfd_{H_{i}}^{G}(X))$$
To see this denote by $\mathcal E$ the Waldhausen category of cofibration sequences in $\Rfd^{G}(X)$, and 
let $\mathcal{E}^{i}$ denote the subcategory of $\mathcal E$ consisting of cofibration sequences 
$$Y'\to Y \to Y/Y'$$ 
such that $Y\in \Rfd_{\leq i}^{G}(X)$, $Y' \in \Rfd_{H_{i}}^{G}(X)$ and $Y/Y' \in \Rfd_{\leq i-1}^{G}(X)$.
By \cite[Propositon 1.3.2]{Wal} the functor $(Y'\to Y \to Y/Y') \to (Y', Y/Y')$ induces a weak equivalence 
$$K(\mathcal{E}^{i}) \overset{\simeq}{\lra} K(\Rfd_{\leq i-1}^{G}(X))\times K(\Rfd_{H_{i}}^{G}(X))$$
 Next, notice that  we have an exact functor $\Phi\colon \mathcal{E}^{i} \to \Rfd_{\leq i}^{G}(X)$ given 
 by $\Phi(Y' \to Y \to Y/Y') = Y$, and an exact functor $\Psi\colon  \Rfd_{\leq i}^{G}(X)\to \mathcal{E}^{i}$ which 
assigns to a space $Y$ the cofibration sequence $\Psi(Y) = (Y'\to Y \to Y/Y')$ where $Y' = X \cup Y^{H_{i}}$.  
These  functors are inverse equivalences of categories and so  they yield a weak equivalence 
$K(\Rfd_{\leq i}^{G}(X))\simeq K(\mathcal{E}^{i})$. 

Since $\Rfd^{G}(X) = \Rfd^{G}_{\leq n}(X)$ arguing by induction we obtain
$$A^{G}(X) = K(\Rfd^{G}(X))  \overset{\simeq}{\lra} \prod_{(H_{i})\in C_{G}} K(\Rfd^{G}_{H_{i}}(X))$$
\end{proof}

Our next goal will be to identify the spectrum $K(\Rfd^{G}_{H}(X))$ for $H\subseteq G$. 
Notice that if $X$ is a $G$-space then the group $WH$ acts the space $X^{H}$. Consider 
the  space $EWH\times X^{H}$ with the action of $WH$ given by $g(y, x) = (gy, yx)$. 
We will show that the following holds: 

\begin{proposition}
\label{RGH EQUIV PROP}
Let $G$ be a finite group, $H\subseteq G$ be a subgroup, and let $X$ be a $G$-space.  
There exists an exact weak equivalence of categories 
$$\Rfd^{G}_{H}(X) \to \Rfd^{WH}(EWH\times X^{H})$$
which induces a weak equivalence of spectra 
$K(\Rfd^{G}_{H}(X)) \overset{\simeq}{\to} A^{WH}(EWH\times X^{H})$. 
\end{proposition}

We will first consider a special case where $H = \{e\}$ is the trivial subgroup of $G$:

\begin{lemma}
\label{RGE EQUIV LEMMA}
Let $G$ be a finite group, and let $X$ be a $G$-space.    
There exists an exact weak equivalence of categories 
$$\Rfd^{G}_{\{e\}}(X) \to \Rfd^{G}(EG\times X)$$
\end{lemma}

\begin{proof}
Consider the exact  functor 
$$\Phi\colon \Rfd_{\{e\}}^{G}(X) \to \Rfd^{G}(EG\times X)$$ 
given by $\Phi(Y) = EG \times Y$ 
for  $Y \in \Rfd_{\{e\}}^{G}(X)$ where $EG\times Y$ is  given the structure of a $G$-retractive space over 
$EG\times X$ in the obvious way. We also have an exact functor 
$$\Psi\colon \Rfd^{G}(EG\times X)  \to \Rfd_{\{e\}}^{G}(X)$$ 
defined as follows. Given a retractive space 
$(r\colon Z \leftrightarrows EG\times X \colon s)\in \Rfd^{G}(EG\times X)$ define
$$\pi_{\ast}Z := \colim(Z \overset{s}{\longleftarrow} EG\times X \overset{\pi}{\lra} X)$$
where $\pi$ is the projection map. The space $\pi_{\ast} Z$ has the structure of a $G$-retractive space over $X$.  
Moreover, since $G$ acts freely on $\pi_{\ast} Z\setminus X$ thus $\pi_{\ast} Z\in \Rfd^{G}_{\{e\}}(X)$. 
We set: $\Psi(Z) = \pi_{\ast}Z$.

For a retractive space  $(r\colon Y \rightleftarrows X \colon s)\in\Rfd_{\{e\}}^{G}(X)$
consider the morphism $\Psi\Phi(Y) \to Y$ in $\Rfd^{G}_{\{e\}}(X)$ induced by
the map of pushouts
\begin{equation*}
\begin{tikzpicture}[baseline=(current bounding box.center)]
\matrix (m) 
[matrix of math nodes, row sep=3em, column sep=4em, text height=1.5ex, text depth=0.25ex]
{
EG\times Y &  EG\times X &  X \\
Y &   X  &  X \\
};
\path[->, thick, font=\scriptsize]
(m-1-1)
edge node[anchor=east] {} (m-2-1)
(m-1-2)
edge node[auto] {$\pi$} (m-1-3)
edge node[anchor=south] {$\id \times s$} (m-1-1)
edge node[anchor=west] {$\pi$}  (m-2-2)
(m-1-3)
edge node[anchor=west] {$=$}  (m-2-3)
(m-2-2) 
edge  node[anchor=north] {$=$} (m-2-3)
edge  node[anchor=north] {$s$} (m-2-1)
; 
\end{tikzpicture}
\end{equation*}
By \cite[II.2.11, Exercise 5]{tomDieck-transfgps} we obtain that  this  morphism is
a weak equivalence in $\Rfd^{G}_{\{e\}}(X)$. 
In this way we get a natural weak equivalence between $\Psi\Phi$ and the identity functor on 
$\Rfd^{G}_{\{e\}}(X)$.  

Finally, for $(r\colon Z \rightleftarrows EG\times X \colon s)\in \Rfd^{G}(EG\times X)$
let  $\varphi_{Z}\colon Z \to \pi_{\ast}Z$ induced by the map of pushouts:
\begin{equation*}
\begin{tikzpicture}[baseline=(current bounding box.center)]
\matrix (m) 
[matrix of math nodes, row sep=3em, column sep=4em, text height=1.5ex, text depth=0.25ex]
{
Z &  EG\times X & EG\times X \\
Z &  EG\times X  &  X \\
};
\path[->, thick, font=\scriptsize]
(m-1-1)
edge node[anchor=east] {$=$} (m-2-1)
(m-1-2)
edge node[auto] {$=$} (m-1-3)
edge node[anchor=south] {$s$} (m-1-1)
edge node[anchor=west] {$=$}  (m-2-2)
(m-1-3)
edge node[anchor=west] {$\pi$}  (m-2-3)
(m-2-2) 
edge  node[anchor=north] {$\pi$} (m-2-3)
edge  node[anchor=north] {$s$} (m-2-1)
; 
\end{tikzpicture}
\end{equation*}
and let $\pi_{EG}\colon EG\times X \to EG$ be the projection map. 
the maps $(\pi_{EG}r, \varphi_{Z})\colon Z \to \Phi\Psi(Z) = EG\times \pi_{\ast}Z$ 
define a natural weak equivalence between the identity functor on 
$\Rfd^{G}(EG\times X)$  and the functor $\Phi\Psi$.
\end{proof}

\begin{proof}[Proof of Proposition \ref{RGH EQUIV PROP}]
By Lemma \ref{RGE EQUIV LEMMA} it will suffice to show that for any subgroup $H\subseteq G$
there exists an exact  equivalence of categories 
$$\Rfd^{G}_{H}(X) \to \Rfd^{WH}_{\{e\}}(X^{H})$$
Let $\Gamma \colon \Rfd_{H}^{G}(X) \to \Rfd_{\{e\}}^{WH}(X^{H})$ be the exact functor 
which associates to a retractive space $Y\in \Rfd_{H}^{G}(X)$ the space $Y^{H}\in \Rfd_{\{e\}}^{WH}(X^{H})$.
We also have an exact functor 
$$\Lambda\colon  \Rfd_{\{e\}}^{WH}(X^{H}) \to \Rfd_{H}^{G}(X)$$
defined as follows. Consider $G/H$ as a left $G$-space and a right $WH$-space. 
Given a retractive space $(r\colon Z \rightleftarrows X^{H} \colon s) \in \Rfd_{\{e\}}^{WH}(X^{H})$
the twisted product $G/H\times_{WH} Z$ is a $G$-retractive space over $G/H\times_{WH} X^{H}$. 
We have a $G$-map 
$$\lambda_{X} \colon G/H\times_{WH} X^{H} \to X$$
given by $\lambda_{X}(gH, x) = gx$. We set:
$$\Lambda(Z) := \colim(G/H\times_{WH}Z \overset{\id\times s} {\longleftarrow} 
G/H\times_{WH} X^{H} \overset{\lambda_{X}}{\lra} X)$$
We claim that $\Gamma$ and $\Lambda$ are inverse equivalences of categories. To see this consider first
a retractive space $(r\colon Y \rightleftarrows X \colon s)\in \Rfd_{H}^{G}(X)$. 
The map of  pushouts diagrams
\begin{equation*}
\begin{tikzpicture}[baseline=(current bounding box.center)]
\matrix (m) 
[matrix of math nodes, row sep=3em, column sep=4em, text height=1.5ex, text depth=0.25ex]
{
G/H\times_{WH}Y^{H} &  G/H\times_{WH} X^{H} & X \\
Y & X  &  X \\
};
\path[->, thick, font=\scriptsize]
(m-1-1)
edge node[anchor=east] {$\lambda_{Y}$} (m-2-1)
(m-1-2)
edge node[auto] {$\lambda_{X}$} (m-1-3)
edge node[anchor=south] {$\id\times s$} (m-1-1)
edge node[anchor=west] {$\lambda_{X}$}  (m-2-2)
(m-1-3)
edge node[anchor=west] {$=$}  (m-2-3)
(m-2-2) 
edge  node[anchor=north] {$=$} (m-2-3)
edge  node[anchor=north] {$s$} (m-2-1)
; 
\end{tikzpicture}
\end{equation*}
induces a map $\eta_{Y}\colon \Lambda\Gamma(Y) \to Y$. 
This map is natural in $Y$, and so defines a natural transformation of functors 
$\eta\colon  \Lambda\Gamma \Ra \rm{Id}_{\Rfd_{H}^{G}(X)}$. We will show that $\eta_{Y}$ is an 
isomorphism for any space $Y$ by constructing its inverse. Take the  $G$-equivariant embedding
$$\mu_{1}\colon s(X) \to \Lambda\Gamma(Y)$$ 
We also have an $NH$-equivariant map $\mu_{2}\colon Y^{H}\to \Lambda\Gamma(Y)$ which is 
the composition of the map $Y^{H}\to G/H\times_{WH}Y^{H}$ that sends $y\in Y^{H}$ to 
the point $[eH, y]$ and the inclusion $G/H\times_{WH}Y^{H} \to \Lambda\Gamma(Y)$. Since 
the sets $s(X)$ and $Y^{H}$ are closed in $Y$, and the maps $\mu_{1}$ and $\mu_{2}$ 
coincide on $s(X)\cap Y^{H} \subseteq Y$ we obtain a continuous map 
$\mu'\colon s(X)\cup Y^{H} \to \Lambda\Gamma(Y)$. Moreover, notice that for every $y\in s(X)\cup Y^{H}$
and for every $g\in G$ such that $gy\in s(X)\cup Y^{H}$ we have $\mu'(gy) = g\mu'(y)$
(we use here the assuption that all orbits of $Y\setminus s(X)$ are isomorphic to $G/H$). 
By \cite[Ch. 1, Theorem 3.3]{Bredon} the map $\mu'$ extends uniquely to a $G$-equivariant map 
$\mu\colon Y \to  \Lambda\Gamma(Y)$. It is straightforward to check that $\mu$ is the inverse of $\eta_{Y}$. 

Construction of a natural isomorphism between the functor $\Gamma\Lambda$ and the identity functor on 
$\Rfd_{\{e\}}^{WH}(X^{H})$ is straightforward. 
\end{proof}

Combining Proposition \ref{AG SPLITTING PROP} with Proposition \ref{RGH EQUIV PROP}
for any $G$-space $X$ we obtain a weak equivalence 
$$ A^{G}(X) \overset{\simeq}{\lra} \prod_{(H)\in C_{G}} A^{WH}(EWH\times X^{H})$$
In order to complete the proof of Theorem \ref{MAIN1 THM} it suffices to show for any subgroup $H\subseteq$
there exists a weak equivalence $A^{WH}(EWH\times X^{H}) \simeq A(EWH\times_{WH}X^{H})$. 
Since the action of $WH$ on $EWH\times X^{H}$ is free this follows from the following fact:

\begin{lemma}
\label{FREE ACTION LEMMA}
If $X$ is a space with a free action of a finite group $G$ then $A^{G}(X)\simeq A(X/G)$. 
\end{lemma}

\begin{proof}
Recall that the spectrum $A(X/G)$ is obtained by applying Waldhausen's $\mathcal{S}_{\bullet}$-construction 
to the category $\Rfd(X/G)$ of homotopy finitely dominated retractive spaces over $X/G$. 
For $(r\colon Y \leftrightarrow X \colon s) \in \Rfd^{G}(X)$ define
$\Gamma(Y) = Y/G$. This space is in a natural way a retractive space over $X/G$
so  it is an object in $\Rfd(X/G)$.  
We also have a functor 
$$\Lambda\colon \Rfd(X/G) \to  \Rfd^{G}(X)$$
defined as follows. Let $q\colon X \to X/G$ be the quotient map. 
For $(r\colon Y \leftrightarrow X/G \colon s) \in \Rfd(X/G)$ let
$$q^{\ast}Y = \{(x, y)\in X\times Y \ | \ q(x) = r(y)\}$$ 
This is a $G$-space with the action of $G$ given by $g(x, y) = (gx, y)$. Moreover, 
$q^{\ast}Y$ is in a natural way a $G$-retractive space over $X$, and it is an object of 
the category $\Rfd^{G}(X)$.  We set: $\Lambda(Y) := q^{\ast}Y$. It is straightforward to check that $\Gamma$
and $\Lambda$ are inverse equivalences of categories. 

\end{proof}


\section{Compatibility with the splitting of  $Q^{G}(X_{+})$}
\label{MODEL OF QGX SEC}

Our next goal is to prove Theorem \ref{MAIN2 THM} which says that the splitting of $A^{G}(X)$ is compatible 
with the splitting of $Q^{G}(X_{+})$, the fixed point spectrum of the equivariant suspension spectrum of $X$. 

Let $G$ be a finite group and let $X$ be a $G$-space. Denote by $\FF^{G}(X)$ the subcategory of 
$\Rfd^{G}(X)$ whose objects are $G$-retractive spaces $r\colon Y \leftrightarrows X \colon  s$ satisfying the 
following conditions:
\benu
\item[(i)] $Y = Y_{X} \sqcup Y_{c}$ such that  $s(X)\subseteq Y_{X}$, 
and the maps $r\colon Y_{X} \leftrightarrows X \colon s$ are inverse $G$-homotopy equivalences. 
\item[(ii)] $Y_{c}$ isomorphic to a finite disjoint union of $G$-spaces of the form $D\times G/H $ where  
$H$ a subgroup of $G$,   $D$ is a contractible space, and 
the $G$-action on $D\times G/H $ is given  by $g(x, kH) = (x, gkH)$.  
\eenu

A morphism between objects $r\colon Y_{X} \sqcup Y_{c} \leftrightarrows X :\!s$
and $r' \colon Y'_{X} \sqcup Y'_{c} \leftrightarrows X :\!s'$ is a $G$-map
$f\colon Y_{X}\sqcup Y_{c}  \to Y_{X}\sqcup Y'_{c}$ 
that satisfies two conditions:
\benu
\item if $y\in Y_{c}$ and $f(y) \in Y_{X}'$ then $f(y) = s'r(y)$;
\item the induced map $f_{\ast}\colon \pi_{0}(f^{-1}(Y_{c}')) \to \pi_{0}(Y'_{c})$ 
is a monomorphism.
\eenu

The category $\FF^{G}(X)$ is a Waldhausen category with weak equivalences and cofibrations inherited from 
$\Rfd^{G}(X)$. The inclusion functor $\FF^{G}(X) \to \Rfd^{G}(X)$ is exact, so it induces a map of spectra
$a^{G}\colon K(\FF^{G}(X)) \to A^{G}(X)$. 

\begin{proposition}
\label{FSPLITTING LEMMA}
There exists a weak equivalence of spectra
$$\tau^{Q}_{X}\colon K(\FF^{G}(X))\overset{\simeq}{\lra} \prod_{(H)\in C_{G}}K(\FF^{WH}(EWH\times X^{H}))$$
natural in $X$  such that the following diagram commutes:
\begin{equation*}
\begin{tikzpicture}[baseline=(current bounding box.center)]
\matrix (m) 
[matrix of math nodes, row sep=3em, column sep=3em, text height=1.5ex, text depth=0.25ex]
{
K(\FF^{G}(X))&   \prod_{(H)\in C_{G}} K(\FF^{WH}(EWH \times X^{H})) &  \\
A^{G}(X) &  \prod_{(H)\in C_{G}} A^{WH}(EWH \times X^{H})\\
};
\path[->, thick, font=\scriptsize]
(m-1-1)
edge node[anchor=east] {$a^{G}$} (m-2-1)
edge node[anchor=north] {$\simeq$} node[anchor=south] {$\tau^{Q}_{X}$} (m-1-2)
(m-2-1)
edge node[anchor=north] {$\simeq$} node[anchor=south] {$\tau^{A}_{X}$} (m-2-2)
(m-1-2)
edge  node[anchor=west] {$\prod a^{WH}$}  (m-2-2)
; 
\end{tikzpicture}
\end{equation*}
\end{proposition}

\begin{proof}
The weak equivalence $\tau^{Q}_{X}$ can be constructed using the same arguments we used in 
Section \ref{PROOF MAIN1 THM SEC} in proofs of Propositions \ref{AG SPLITTING PROP} and
\ref{RGH EQUIV PROP}, taking the category $\FF^{G}(X)$ in place of $\Rfd^{G}(X)$ and 
the category $\FF^{G}_{H}(X) = \FF^{G}(X)\cap \Rfd^{G}_{H}(X)$ in place of $\Rfd^{G}_{H}(X)$. 
Commutativity of the diagram follows directly from this construction. 
\end{proof}

In order to complete the proof of Theorem \ref{MAIN2 THM} it will suffice to show the following holds: 
\benu
\item[(i)] $K(\FF^{G}(X))\simeq Q^{G}(X_{+})$ 
\item[(ii)] $K(\FF^{WH}(EWH \times X^{H})) \simeq Q(EWH\times_{WH}X^{H}_{+})$ for each $H\in G$.
\eenu
Notice that it is enough to prove the second of these statements. Indeed, it will give
$K(\FF^{G}(X))\simeq \prod_{(H)\in \FF_{G}} Q(EWH\times_{WH}X^{H}_{+})$,
and by \cite{tomDieck-Orbitttypen}, \cite[Theorem 11.1]{Lewis-May}.
the product of on the right hand side is equivalent to $Q^{G}(X_{+})$

The proof of (ii) will be split in two steps. First, for a space $X$  let  $\FF(X)$ denote the category whose objects 
are (non-equivariant) retractive space $r\colon  Y \leftrightarrows X :\!s$ satisfying the following conditions: 
\benu
\item $Y = Y_{X} \sqcup \bigsqcup_{i\in I} B_{i}$ where $I$ is a finite set and $B_{i}$ is a contractible space for each $i\in I$; 
\item the map $s\colon X \to Y$ is a cofibration and  $s(X) \subseteq  Y_{X}$; 
\item the maps $s\colon X\to Y_{X}$ and $r\colon Y_{X} \to X$ are inverse homotopy equivalences. 
\eenu

A morphism between objects $r\colon Y_{X} \sqcup \bigsqcup_{i\in I} B_{i} \leftrightarrows X :\!s$
and $r' \colon Y_{X} \sqcup \bigsqcup_{i\in I} B_{i} \leftrightarrows X :\!s'$ is a continuous map 
$f\colon Y_{X}\sqcup \bigsqcup_{i\in I} B_{i} \to Y'_{X}\sqcup \bigsqcup_{j\in J} B'_{j}$ over and under $X$ 
that satisfies two conditions:
\benu
\item if $f(B_{i})\subseteq Y_{X}'$ then $f|_{B_{i}} = s'r$. 
\item the induced map $f_{\ast}\colon \pi_{0}(f^{-1}(\bigsqcup_{j\in J} B'_{j})) \to \pi_{0}(\bigsqcup_{j\in J} B'_{j})$ 
is a monomorphism.
\eenu

The category $\FF(X)$ is a Waldhausen category where a morphism $f\colon Y\to Y'$ is a weak equivalence 
if it is a weak homotopy equivalences and $f$ is a cofibration it has the homotopy 
extension property.

\begin{lemma}
\label{FGE EQUIV LEMMA}
Let $G$ be a finite group, and let $X$ be a $G$-space. If $G$ acts freely on $X$ then there exists a weak equivalence 
of spectra
$K(\FF^{G}(X)) \simeq  K(\FF(X/G))$
\end{lemma}

The proof is essentially the same as the proof of Lemma \ref{FREE ACTION LEMMA}.

Going back to claim (ii), since the action of $WH$ on $EWH\times X^{H}$ is free by Lemma \ref{FGE EQUIV LEMMA}
we obtain a weak equivalence  $K(\FF^{WH}(EWH \times X^{H})) \simeq K(\FF(EWH\times_{WH} X^{H}))$. 
In order to complete the proof of existence of the weak equivalence  (ii), and thus complete the proof
of Theorem \ref{MAIN2 THM} in remains to show that the following holds:

\begin{lemma}
For any space $X$ there is a weak equivalence of spectra
$$K(\FF(X)) \simeq Q(X_{+})$$
\end{lemma}

\begin{proof}
We will assume that $X$ is connected space. The general case follows from essentially the same argument. 
Notice that for any cofibration sequence $Y'\to Y \to Y/Y'$ in $\FF(X)$ there is a functorial weak equivalence 
$Y \simeq Y'\sqcup Y/Y'$. This implies that we have a weak equivalence of spectra
$$K(\FF(X))\simeq \Omega |w\mathcal{N}_{\bullet}\FF(X)|$$
where $w\mathcal{N}_{\bullet}\FF(X)$ is a simplicial category defined as in \cite[\S 1.8]{Wal}. 
Let $L_{X}$ denote the loop group of $X$. Let $\mathcal{R}^{0}_{k}(\ast, L_{X})$ be the category of 
$0$-spherical objects of rank $k$ defined as in \cite[p. 386]{Wal} and let 
$\mathcal{R}^{0}(\ast, L_{X}) = \bigcup_{k} \mathcal{R}^{0}_{k}(\ast, L_{X})$. We have:
$$\Omega |w\mathcal{N}_{\bullet}\FF(X)| \simeq \Omega|w\mathcal{R}^{0}(\ast, L_{X})|
\simeq \Z\times \lim_{k}|w\mathcal{R}^{0}_{k}(\ast, L_{X})|^{+}\simeq Q(X_{+})$$
The first of these weak equivalences can be obtained by essentially the same argument as one 
used in the proof of \cite[Proposition 2.1.4]{Wal}, while the other two comes from a theorem of 
Segal \cite{Segal-Gamma} (see also \cite[p.386]{Wal}). 
\end{proof}


\bibliographystyle{plain}
\bibliography{splitting}

\end{document}